\documentclass{amsart}

\usepackage{amsfonts, amsmath, amssymb, amscd}
\usepackage{pstricks,pst-node,pst-plot}
\usepackage{epic}
\usepackage{autograph}

\newtheorem{theorem}{Theorem}[section]

\newtheorem{proposition}[theorem]{Proposition}
\newtheorem{corollary}[theorem]{Corollary}

\theoremstyle{definition}

\theoremstyle{plain}

\begin{document}
    \title{Rational Subsets of Polycyclic Monoids and Valence Automata}

\maketitle
\begin{center}
    Elaine Render and Mark Kambites \\

    \medskip

    School of Mathematics, \ University of Manchester \\
    Manchester M60 1QD, \  England \\

    \medskip

    \texttt{E.Render@maths.manchester.ac.uk} \\
    \texttt{Mark.Kambites@manchester.ac.uk} \\

    \medskip

\end{center}

\begin{abstract}
We study the classes of languages defined by valence automata with
rational target sets (or equivalently, regular valence grammars with
rational target sets), where the valence monoid is drawn from the
important class of polycyclic monoids. We show that for polycyclic
monoids of rank $2$ or more, such automata accept exactly the
context-free languages. For the polycyclic monoid of rank $1$ (that
is, the bicyclic monoid), they accept a class of languages strictly
including the partially blind one-counter languages. Key to the
proof is a description of the rational subsets of polycyclic and
bicyclic monoids, other consequences of which include the
decidability of the rational subset membership problem, and the
closure of the class of rational subsets under intersection and
complement.
\end{abstract}

\bigskip

\section{Introduction}

Both mathematicians and computer scientists have found applications
for finite automata augmented with registers which store values from
a given group or monoid, and are modified by multiplication. These
automata, variously known as \textit{valence automata},
\textit{extended finite automata} or \textit{$M$-automata}, provide
an algebraic method to characterize important language classes such
as the context-free, recursively enumerable and blind counter
languages (see \cite{Gilman96}). Their study provides insight into
computational problems in algebra (see, for example,
\cite{Kambites/Silva/Steinberg07}). These automata are also closely
related to \textit{regulated rewriting systems}, and in particular
the \emph{valence grammars} introduced by Paun \cite{Paun80}: the
languages accepted by $M$-automata are exactly the languages
generated by regular $M$-valence grammars \cite{Fernau/Stiebe02}.

Traditionally, the monoid registers are initialised to the identity
element, and a word is accepted only if it can be read by a
successful computation which results in the register being returned
to the identity. Several authors have observed that the power of
these automata to describe language classes may be increased by
allowing a more general set of accepting values in the register.
Fernau and Stiebe \cite{Fernau/Stiebe01} began the systematic study
of the resulting \textit{valence automata with target sets}, along
with the corresponding class of regulated grammars. In particular
they considered the natural restriction that the target set be a
\textit{rational subset} of the register monoid.

Of particular interest, when considering semigroups and monoids in
relation to automata theory, is the class of \textit{polycyclic
monoids}. The polycyclic monoid of rank $n$ is the natural algebraic
model of a pushdown store on an $n$-letter alphabet. For $M$ a
polycyclic monoid of rank $2$ or more, it is well-known that
$M$-automata are equivalent to pushdown automata, and hence accept
exactly the context-free languages. The polycyclic monoid of rank
$1$ is called the \textit{bicyclic monoid} and usually denoted $B$;
we shall see below that $B$-automata accept exactly the
\textit{partially blind one-counter languages} defined by Greibach
\cite{Greibach78}.

One of the main objectives of this paper is to consider the class of
languages accepted by polycyclic monoid valence automata with rational
target sets. It transpires that, for $M$ a polycyclic monoid of rank $2$
or more, every language accepted by an $M$-automaton with rational target
set is context-free, and hence is accepted by an $M$-automaton with
target set $\lbrace 1 \rbrace$. In the rank $1$ case the situation is rather
different; a language accepted by a $B$-automata with rational target set
need not be a partially blind one-counter language, but it is always a
finite union of languages, each of which is the concatenation of two
partially blind one-counter languages.

A key element of the proofs is a simple but extremely useful characterisation
of the rational subsets of polycyclic monoids (Corollary~\ref{splitcor} below).
From this we are easily able to derive a number of other consequences which may be of
independent interest. These include the facts that the rational subsets
of a finitely generated polycyclic monoid form a boolean algebra (with
operations effectively computable), and that membership is uniformly
decidable for rational subsets of polycyclic monoids.

In addition to this introduction, the present paper is divided into four
sections. Section~\ref{preliminaries} recalls some basic definitions from
formal language theory and the theory of valence automata, while
Section~\ref{transductions} establishes some foundational results about
valence automata with rational target sets.
In Section~\ref{adjoinzerosection} we consider the effect of adjoining a
zero to a monoid $M$ upon the classes of languages accepted by $M$-automata
and by $M$-automata with rational target sets. Finally, in
Section~\ref{polycyclicsection} we turn our attention to polycyclic and
bicyclic monoids, proving our main results about both rational subsets
and valence automata with rational target sets.

\section{Preliminaries}\label{preliminaries}

Firstly, we recall some basic ideas from formal language theory. Let
$\Sigma$ be a finite alphabet. Then we denote by $\Sigma^*$ the set
of all words over $\Sigma$ and by $\epsilon$ the empty word. Under
the operation of concatenation and with the neutral element
$\epsilon$, $\Sigma^*$ forms a free monoid. A \emph{finite
automaton} over $\Sigma^*$ is a finite directed graph with each edge
labelled with an element of $\Sigma^*$, and with a distinguished
initial vertex and a set of distinguished terminal vertices. A word
$w \in \Sigma^*$ is accepted by the automaton if there exists some
path connecting the initial vertex with some terminal vertex, the product
of whose edge labels in order is $w$. The set of all words accepted by the
automaton is denoted $L$ or for an automaton $A$ sometimes $L(A)$,
and is called the \emph{language} accepted by $A$. A language
accepted by a finite automaton is called \emph{rational} or
\emph{regular}.

More generally, if $M$ is a monoid then a \textit{finite automaton
over $M$} is a finite directed graph with each edge labelled with an
element of $M$, and with a distinguished initial vertex and a set of
distinguished terminal vertices. An element $m \in M$ is accepted by
the automaton if there exists some path connecting the initial
vertex with some terminal vertex, the product in order of whose edge
labels is $m$. The \emph{subset accepted} is the set of all elements
accepted; a subset of $M$ which is accepted by some finite automaton
is called a \emph{rational subset}. The rational subsets of $M$ are
exactly the homomorphic images in $M$ of regular languages.

We now recall the definition of a finite valence automaton, or
$M$-automaton. Let $M$ be a monoid with identity $1$ and let
$\Sigma$ be an alphabet. An \emph{$M$-valence automaton} (or
$M$-automaton for short) over $\Sigma$ is a finite automaton over
the direct product $M \times \Sigma^*$. We say that it accepts a
word $w \in \Sigma^*$ if it accepts $(1, w)$, that is if there
exists a path connecting the initial vertex to some terminal vertex
labelled $(1,w)$.

Intuitively, we visualize an $M$-automaton as a finite automaton
augmented with a memory register which can store an element of $M$;
the register is initialized to the identity element, is modified by
right multiplication by elements of $M$, and for a word to be
accepted the element present in the memory register on completion
must be the identity element. We write $F_1(M)$ for the class of all
languages accepted by $M$-automata, or equivalently for the class of
languages generated by $M$-valence grammars \cite{Fernau/Stiebe02}.
More generally, an \emph{$M$-automaton with (rational) target set} is an
$M$-valence automaton together with a (rational) subset $X \subseteq
M$. A word $w \in \Sigma^*$ is accepted by such an automaton if it
accepts $(x,w)$ for some $x \in X$. We denote by $F_{Rat}(M)$ the
family of languages accepted by $M$-automata with rational target
sets. We recall the following result of Fernau and Stiebe
\cite{Fernau/Stiebe01}.
\begin{theorem}[Fernau and Stiebe 2001] \label{GXX}
Let $G$ be a group. Then $F_{Rat}(G) = F_{1}(G)$.
\end{theorem}

\section{Automata, Transductions and Closure Properties}\label{transductions}

In this section we study the relationship between rational
transductions and $M$-automata with target sets. Consider a finite
automaton over the direct product $\Omega^* \times \Sigma^*$. We
call an automaton of this type a \emph{rational transducer} from
$\Omega$ to $\Sigma$; it recognises a relation $R \subseteq \Omega^*
\times \Sigma^*$ called a \emph{rational transduction}. The image of
a language $L \subseteq \Omega^*$ under the relation $R$ is the set
of $y \in \Sigma^*$ such that $(x,y) \in R$ for some $x \in L$. We
say that a language $K$ is a \emph{rational transduction of} a
language $L$ if $K$ is the image of $L$ under some rational
transduction. The following is a straightforward generalisation of a
well-known observation concerning $M$-automata (see for example
\cite[Proposition~2]{Kambites06}).

\begin{proposition} \label{transduction}
Let $X$ be a subset of a monoid $M$, and let $L \subseteq \Sigma^*$
be a regular language. Then the following are equivalent:
\begin{itemize}
\item[(i)] $L$ is accepted by an $M$-automaton with target set $X$;
\item[(ii)] there exists a finite alphabet $\Omega$ and a morphism
$\omega : \Omega^* \to M$ such that $L$ is a rational transduction of
$X \omega^{-1}$.
\end{itemize}
If $M$ is finitely generated then the following condition is also
equivalent to those above.
\begin{itemize}
\item[(iii)] for every finite choice of generators $\omega : \Omega^* \to M$
for $M$, $L$ is a rational transduction of $X \omega^{-1}$.
\end{itemize}
\end{proposition}
\begin{proof}
To show that (i) implies (ii), suppose $L$ is accepted by an
$M$-automaton with target set $X$. Choose a finite alphabet $\Omega$
and a map $\omega: \Omega^* \to M$ such that the image $\Omega^*
\omega$ contains every element of $M$ which forms the first
component of an edge-label in the automaton. We now obtain from the
automaton a transducer from $\Omega$ to $\Sigma$ by replacing each
edge label $(m, x)$ with $(w, x)$ where $w \in \Omega^*$ is some
word such that $w \omega = m$. It is a routine exercise to verify
that $L$ is the image of $X\omega^{-1}$ under the given
transduction.

Conversely, suppose we are given a map $\omega : \Omega^* \to M$ and
a transducer from $\Omega$ to $\Sigma$. We construct from the
transducer an $M$-automaton with target set $X$ by replacing each
edge label of the form $(w, x)$ with $(w \omega, x)$. It is readily
verified that the language accepted by this $M$-automaton is exactly
the image of $X$ under the transduction.

Suppose now that $M$ is finitely generated. Clearly, (iii) implies
(ii). Finally, if (ii) holds then we can extend $\omega$ arbitrarily
to a finite choice of generators $\omega' : (\Omega')^* \to M$ for
$M$, and check that we still have the desired property, so that
(iii) holds.
\end{proof}

In particular, Proposition~\ref{transduction} gives a characterisation
in terms of rational subsets and transductions of each class of languages
accepted by $M$-automata with rational target sets.

\begin{proposition} \label{Mrat}
Let $M$ be a monoid and $L \subseteq \Sigma^*$ a language. Then the
following are equivalent.
\begin{itemize}
\item[(i)] $L \in F_{Rat}(M)$;
\item[(ii)] there exists a finite alphabet $\Omega$, a morphism
$\omega : \Omega^* \to M$ and a rational subset $X \subseteq M$ such
that $L$ is a rational transduction of $X \omega^{-1}$.
\end{itemize}
If $M$ is finitely generated then the following condition is also
equivalent to those above.
\begin{itemize}
\item[(iii)] there exists a rational subset $X \subseteq M$ such that
for every finite choice of generators $\omega : \Omega^* \to M$ for
$M$, $L$ is a rational transduction of $X \omega^{-1}$.
\end{itemize}
\end{proposition}

Recall that a \textit{rational cone} (also known as a \textit{full
trio}) is a family of languages closed under rational transduction,
or equivalently under morphism, inverse morphism, and intersection
with regular languages \cite[Section~V.2]{Berstel79}. Since rational
transductions are closed under composition
\cite[Theorem~III.4.4]{Berstel79} we have the following immediate
corollary.
\begin{corollary}\label{ratcone}
$F_{Rat}(M)$ is a rational cone. In particular, it is closed under
morphism, inverse morphism, intersection with regular languages, and
(since it contains a non-empty language) union with regular
languages.
\end{corollary}

\section{Adjoining a Zero}\label{adjoinzerosection}

In this section we show that adjoining a zero to a monoid $M$ makes
no difference to the families of languages accepted either by
$M$-automata or by $M$-automata with rational target sets. Recall that if
$M$ is a monoid, the result of \textit{adjoining a zero} to $M$ is
the monoid $M^0$ with set of elements $M \cup \lbrace 0 \rbrace$
where $0$ is a new symbol not in $M$, and multiplication given by
$$st = \begin{cases} \text{ the $M$-product } st &\text{ if } s, t \in M \\
                     0 &\text{ otherwise.}
       \end{cases}$$
We begin with the $M$-automaton case, where the required result is a very
simple observation.
\begin{proposition}\label{M0M}
Let $M$ be a monoid. Then $F_1(M^0) = F_1(M)$.
\end{proposition}
\begin{proof}
That $F_1(M) \subseteq F_1(M^0)$ is immediate, so we need only prove
the converse. Suppose $L \in F_1(M^0)$, and let $A$ be an
$M^0$-automaton accepting $L$. Clearly any path in $A$ containing an
edge with first label component $0$ will itself have first label
component $0$; thus, no accepting path in $A$ can contain such an
edge. It follows that by removing all edges whose label has first
component $0$, we obtain a new $M^0$-automaton $B$ accepting the
language $L$. But now since $M$ is a submonoid of $M^0$, $B$ can be
interpreted as an $M$-automaton accepting $L$, so that $L \in
F_1(M)$ as required.
\end{proof}

Next we establish the corresponding result for $M$-automata with
rational target sets, which is a little more involved.

\begin{theorem} \label{G0G}
Let $M$ be a monoid. Then $F_{Rat}(M^0) = F_{Rat}(M)$.
\end{theorem}
\begin{proof}
That $F_{Rat}(M) \subseteq F_{Rat}(M^0)$ is immediate. For the
converse, suppose $L \in F_{Rat}(M^0)$. Then we may choose an
$M^0$-automaton $A$ accepting $L$ with rational target set $X \subseteq M$.

Let $L_0$ be the language of words $w \subseteq \Sigma^*$ such that
$(0, w)$ labels a path from the initial vertex to a terminal vertex.
Let $L_1$ be the set of words $w$ such that $(m, w)$ labels a path
from the initial vertex to a terminal vertex for some $m \in X
\setminus \lbrace 0 \rbrace$. Clearly either $L = L_0 \cup L_1$ (in
the case that $0 \in X$) or $L = L_1$ (if $0 \notin X$). We claim
that $L_0$ is regular and $L_1 \in F_{Rat}(M)$. By
Proposition~\ref{ratcone} this will suffice to complete the proof.

The argument to show that $L_1 \in F_{Rat}(M)$ is very similar to
the proof of Proposition~\ref{M0M}. We construct from the
$M^0$-automaton $A$ a new $M$-automaton $B$ by simply removing each
edge with label of the form $(0, m)$. The new automaton $B$ has
target set $X \setminus \lbrace 0 \rbrace$. It is straightforward to
show that $B$ accepts exactly the language $L_1$.

It remains to show that $L_0$ is regular. Let $Q$ be the vertex set
of the automaton $A$, and let $Q_{0} = \lbrace q_0 \mid q \in Q \rbrace$ and
$Q_1 = \lbrace q_1 \mid q \in Q \rbrace$ be disjoint copies of $Q$. We
define from $A$ a finite automaton $C$ with
\begin{itemize}
\item vertex set $Q_0 \cup Q_1$;
\item for each edge in $A$ from $p$ to $q$ with label of the form $(m, x)$
\begin{itemize}
\item an edge from $p_0$ to $q_0$ labelled $x$ and
\item an edge from $p_1$ to $q_1$ labelled $x$;
\end{itemize}
\item for each edge in $A$ from $p$ to $q$ with label of the form $(0, x)$
\begin{itemize}
\item an edge from $p_0$ to $q_1$ labelled $x$ and
\item an edge from $p_{1}$ to $q_{1}$ labelled $x$;
\end{itemize}
\item initial vertex $q_0$ where $q$ is the initial vertex of $A$; and
\item terminal vertices $q_1$ whenever $q$ is a terminal vertex of $A$.
\end{itemize}
We shall show that $C$ accepts exactly the language $L_{0}$. Let $w \in
L_{0}$. Then there exists an accepting path $\pi$ through $A$ labelled
$(0,w)$. It follows from the definition of $M^0$ that no product of
non-zero elements can equal $0$; hence, this path must traverse at least
one edge with label of the form $(0, x)$ for some $x \in \Sigma^*$.
 Suppose then that $\pi = \pi_1 \pi_2 \pi_3$ where $\pi_1$ is a path
from the initial vertex to a vertex $p$ with label $(m_1, w_1)$,
$\pi_2$ is an edge from $p$ to a vertex $q$ with label $(0, x)$, and
$\pi_3$ is a path from $q$ to a terminal vertex with label $(m_3, w_3)$.
It follows easily from the definition of $C$ that it has a path from the
initial vertex to $p_0$ labelled $w_1$, an edge from $p_0$ to $q_1$ with
label $x$, and an edge from $q_1$ to a terminal vertex with label $w_3$.
Hence, $w = w_1 x w_3$ is accepted by $C$, as required.

Conversely suppose $w \in L(C)$, and let $\pi$ be an accepting path
for $w$. Notice that the initial vertex of $C$ lies in $Q_0$ while
all the terminal vertices lie in $Q_1$. Then $\pi = \pi_1 \pi_2
\pi_3$ where $\pi_1$ is a path from the initial vertex to some $p_0$
with label $w_1$, $\pi_2$ is an edge from $p_0$ to some $q_1$ with
label $x$, $\pi_3$ is a path from $q_1$ to a terminal vertex with
label $w_3$ where $w = w_1 x w_3$. Now it follows easily from the
definition of $C$ that $A$ has paths from the initial vertex to $p$
with label of the form $(m_1, w_1)$, from $p$ to $q$ with label $(0,
x)$ and from $q$ to a terminal vertex with label of the form $(m_3,
w_3)$. Thus, $A$ accepts $(m_1 0 m_3, w_1 x w_3) = (0, w)$ so that
$w \in L_0$ as required.
\end{proof}

Combining Theorems~\ref{G0G} with the result of Fernau and Stiebe
\cite{Fernau/Stiebe01} mentioned above (Theorem~\ref{GXX}) gives us
the following immediate corollary.
\begin{corollary}
Let $G$ be a group. Then
$$F_{Rat}(G^0) = F_{Rat}(G) = F_1(G) = F_1(G^0).$$
\end{corollary}

\section{Polycyclic Monoids}\label{polycyclicsection}

In this section we study the language classes $F_1(M)$ and
$F_{Rat}(M)$, where $M$ is drawn from the class of
\textit{polycyclic monoids}, which form the natural algebraic models
of pushdown stores. In the process, we obtain a number of
results about rational subsets of these monoids which may be of
independent interest.

Let $X$ be a set. Recall that the \emph{polycyclic monoid} on $X$ is
the monoid $P(X)$ generated, under the operation of relational
composition, by the partial bijections of the form
\[ p_{x}: X^* \to X^*, \quad w \mapsto wx\]
and
\[ q_{x}: X^*x \to X^*, \quad wx \mapsto w. \]
The monoid $P(X)$ is a natural algebraic model of a pushdown store
on the alphabet $X$, with $p_x$ and $q_x$ corresponding to the
elementary operations of \emph{pushing $x$} and \emph{popping $x$}
(where defined) respectively, and composition to performing these
operations in sequence. For a more detailed introduction see
\cite{Kambites06}.

Clearly for any $x \in X$, the composition $p_x q_x$ is the identity
map. On the other hand, if $x$ and $y$ are distinct letters in $X$,
then $p_x q_y$ is the \textit{empty map} which constitutes a zero
element in $P(X)$. In the case $|X| = 1$, say $X = \lbrace x
\rbrace$, the monoid $P(X)$ is called the \textit{bicyclic monoid},
and is often denoted $B$. The partial bijections $p_x$ and $q_x$
alone (which we shall often denote just $p$ and $q$) do not generate
the empty map, and so the bicyclic monoid does not have a zero
element; to avoid having to treat it as a special case, it is
convenient to write $P^0(X)$ for the union of $P(X)$ with the empty
map; thus we have $P^0(X) = P(X)$ if $|X| \geq 2$ but $P^0(X)$
isomorphic to $P(X)$ with a zero adjoined if $|X| = 1$.

Let $P_X = \lbrace p_x \mid x \in X \rbrace$ and $Q_X = \lbrace q_x
\mid x \in X \rbrace$, and let $z$ be a new symbol not in $P_X \cup
Q_X$ which will represent the zero element. Let $\Sigma_X = P_X \cup
Q_X \cup \lbrace z \rbrace$. Then there is an obvious surjective
morphism $\sigma : \Sigma_X^* \to P^0(X)$, and indeed $P^0(X)$
admits the monoid presentation
\begin{align*}
P^0(X) = \langle & \Sigma_X \mid p_{x}q_{x} = 1, p_{x}q_{y}=z, \\
& zp_x = z q_x = p_x z = q_x z = zz = z \ \textrm{for all}\ x,y \in
X,\ x\neq y  \rangle.
\end{align*}
It is well-known (see for example
\cite{Gilman96,Kambites06}) that for $|X| \geq 2$, a
$P(X)$-automaton is equivalent to a pushdown automaton with stack
alphabet $X$, so that the language class $F_1(P(X))$ is exactly the
class of context-free languages. Greibach \cite{Greibach78} has
introduced and studied the class of \textit{partially blind counter
automata}. The latter are non-determinstic finite automata augmented
with a number of non-negative integer counters which can be
incremented and decremented but not read; attempting to decrement a
counter whose value is $0$ causes the computation to fail. The
counters are initialized to 0, and a word is accepted only if some
computation reading that word places the finite state control in an
accepting state and returns all counters to 0. The following
equivalence follows immediately from the definitions.
\begin{proposition}\label{bicyclic}
For any $n > 0$, $F_{1}(B^n)$ is exactly the class of languages
accepted by partially blind $n$-counter automata.
\end{proposition}

We now turn our attention to the classes $F_{Rat}(P(X))$ of
languages accepted by polycyclic monoid automata with rational
target sets. For $|X| \geq 2$, it transpires that every language
accepted by a $P(X)$-automaton with rational target set is accepted
by a $P(X)$-automaton, and hence that $F_{Rat}(P(X))$ is the class
of context free languages. In order to prove this, we shall need
some results about rational subsets of polycyclic monoids, which we
establish using techniques from string rewriting theory.

Recall that a  \emph{monadic rewriting system} $\Lambda$ over an
alphabet $\Sigma$ is a subset of $\Sigma^* \times \{\Sigma \cup
\{\epsilon\}\}$. We normally write an element $(w,x) \in \Lambda$ as
$w \to x$. Then we write $u \Rightarrow v$ if $u=rws \in \Sigma^*$
and $v = rxs \in \Sigma^*$ with $w \to x$. Denote by $\Rightarrow^*$
the transitive, reflexive closure of the relation $\Rightarrow$. If
$u \Rightarrow^* v$ we say that $u$ is an \emph{ancestor} of $v$
under $\Lambda$ and $v$ is a \emph{descendant} of $u$ under
$\Lambda$; we write $L \Lambda$ for the set of all descendants of
words in $L$. It is well-known that if $L$ is regular then $L
\Lambda$ is again a regular language; if moreover the rewriting
system $\Lambda$ is finite, a finite automaton recognising $L
\Lambda$ can be effectively computed from a finite automaton
recognising $L$. For more information on such systems see
\cite{Book/Jantzen/Wrathall82,Book/Otto93}.

\begin{theorem}\label{q*p*}
Let $X$ be a finite alphabet and $R$ a rational subset of $P^0(X)$.
Then there exists a regular language
$$L \subseteq Q_X^*P_X^* \cup \{z\}$$
such that $L \sigma = R$. Moreover, there is an algorithm which,
given an automaton recognizing a regular language $G \subseteq
\Sigma_X^*$, constructs an automaton recognising a language $L
\subseteq Q_X^*P_X^* \cup \{z\}$ with $L \sigma = G \sigma$.
\end{theorem}
\begin{proof}
Since $R$ is rational, there exists a regular language $K \subseteq
\Sigma_X^*$ such that $K\sigma = R$. We define a monadic rewriting
system $\Lambda$ on $\Sigma_X^*$ with the following rules:
\begin{eqnarray*}
p_{x}q_{x} \rightarrow \epsilon, \quad &p_{x}q_{y} \rightarrow z,& \quad zq_{x} \rightarrow z,\\
p_{x}z \rightarrow z, \quad &zp_{x}
\rightarrow z,& \quad q_{x}z \rightarrow z,  \\
&zz \rightarrow z &
\end{eqnarray*}
for all $x,y \in X$ with $x \neq y$.
Notice that the language of $\Lambda$-irreducible words is
exactly $Q_X^* P_X^* \cup \lbrace z \rbrace$. With this in mind, we
define
$$L \ = \ K \Lambda \cap (Q_X^* P_X^* \cup \lbrace z \rbrace)$$
Certainly $L$ is regular, and moreover an automaton for $L$ can be
effectively computed from an automaton for $K$. Thus, it will
suffice to show that $L \sigma = R$.

By definition $L \sigma \subseteq (K\Lambda)\sigma$, and since the
rewriting rules are all relations satisfied in $P^0(X)$,
$$(K \Lambda) \sigma \ \subseteq \ K \sigma \ = \ R.$$
Conversely, if $s \in R$ then $s = w \sigma$ for some $w \in K$. Now
the rules of $\Lambda$ are all length-reducing, so $w$ must clearly
have an irreducible descendant, say $w'$. But now $w' \in L$ and $w' \sigma
= w \sigma = s$ so that $s \in L \sigma$. Thus, $L \sigma = R$ as
required.
\end{proof}

As an immediate corollary, we obtain a corresponding result for
bicyclic monoids
\begin{corollary}\label{bicyclicq*p*}
Let $R$ be a rational subset of a bicyclic monoid $B$, and $\sigma :
\lbrace p, q \rbrace^* \to B$ the natural morphism. Then there
exists a regular language $L \subseteq q^*p^*$ such that $L \sigma =
R$. Moreover, there is an algorithm which, given an automaton
recognizing a regular language $G \subseteq \{ p, q \}^*$,
constructs an automaton recognising a language $L \subseteq q^*p^*$
with $L \sigma = G \sigma$.
\end{corollary}

Before proceeding to apply the theorem to polycyclic monoid automata
with target sets, we note some general consequences of
Theorem~\ref{q*p*} for rational subsets of polycyclic monoids.
Recall that a collection of subsets of a given \textit{base set} is
called a \textit{boolean algebra} if it is closed under union,
intersection and complement within the base set.
\begin{corollary}\label{booleanalgebra}
The rational subsets of any finitely generated polycyclic monoid form a boolean
algebra. Moreover, the operations of union, intersection and complement are
effectively computable.
\end{corollary}
\begin{proof}
The set of rational subsets of a monoid is always (effectively)
closed under union, as a simple consequence of non-determinism.
Since intersection can be described in terms of union and
complement, it suffices to show that the rational subsets of
polycyclic monoids are closed (effectively) under complement. To
this end, suppose first that $R$ is a rational subset of a finitely generated
polycyclic monoid $P(X)$ with $|X| \geq 2$. Then by Theorem~\ref{q*p*}, there is a
regular language $L \subseteq (Q_X^* P_X^* \cup \{z\})$  such that
$L \sigma = R$. Let $K = (Q_X^* P_X^* \cup \{z\}) \setminus L$. Then
$K$ is regular and, since $Q_X^* P_X^* \cup \{z\}$ contains a
unique representative for every element of $P(X)$, it is readily
verified that $K \sigma = P(X) \setminus (L \sigma)$. Thus, $P(X)
\setminus (L \sigma)$ is a rational subset of $P(X)$, as required.

For effective computation of complements, observe that given an
automaton recognizing a language $R = \Sigma_X^*$, we can by
Theorem~\ref{q*p*} construct an automaton recognizing a regular
language $L \subseteq (Q_X^* P_X^* \cup \{z\})$ with $L \sigma = R
\sigma$. Clearly we can then compute the complement $K = (Q_X^*
P_X^* \cup \{z\}) \setminus L$ of $L$ in $(Q_X^* P_X^* \cup \{z\})$,
and since $K \sigma = P(X) \setminus (L \sigma)$, this suffices.

In the case that $|X| = 1$, the statement can be proved in a similar way
but using Corollary~\ref{bicyclicq*p*} in place of
Theorem~\ref{q*p*}.
\end{proof}

Recall that the \textit{rational subset problem} for a monoid $M$ is
the algorithmic problem of deciding, given a rational subset of $M$
(specified using an automaton over a fixed generating alphabet) and
an element of $M$ (specified as a word over the same generating
alphabet), whether the given element belongs to the given subset.
The decidability of this problem is well-known to be independent of
the chosen generating set \cite[Corollary~3.4]{Kambites/Silva/Steinberg07}. As
another corollary, we obtain the decidability of this problem for
finitely generated polycyclic monoids.
\begin{corollary}
Finitely generated polycyclic monoids have decidable rational subset
problem.
\end{corollary}
\begin{proof}
Let $|X| \geq 2$ [respectively, $|X| = 1$]. Suppose we are given a
rational subset $R$ of $P(X)$ (specified as an automaton over
$\Sigma_X^*$ [respectively $\lbrace p, q \rbrace^*$]) and an element
$w$ (specified as a word in the appropriate alphabet). Clearly, we
can compute $\lbrace w \rbrace$ as a regular language. Now by
Corollary~\ref{booleanalgebra} we can compute a regular language $K
\subseteq \Sigma_X^*$ [respectively, $\lbrace p, q \rbrace^*$] such
that $K \sigma = R \cap \lbrace w \rbrace \sigma$. So $w \sigma \in
R$ if and only if $R \cap \lbrace w \rbrace$ is non-empty, that is,
if and only if $K$ is non-empty. Since emptiness of regular
languages is testable, this completes the proof.
\end{proof}

We now return to our main task of proving that $F_{Rat}(M) = F_1(M)$
for $M$ a polycyclic monoid of rank $2$ or more, that is, that
polycyclic monoid automata with target sets accept only context-free
languages. We shall need some preliminary results.

\begin{corollary}\label{splitcor}
Let $R$ be a rational subset of $P^0(X)$ and suppose that $0 \notin
R$. Then there exists an integer $n$ and regular languages $Q_1,
\dots, Q_n \subseteq Q_X^*$ and $P_1, \dots, P_n \subseteq P_X^*$
such that
$$R = \bigcup_{i=1}^n \left( Q_i P_i\right) \sigma.$$
\end{corollary}

\begin{proof}
By Theorem~\ref{q*p*}, there is a regular language $L \subseteq
Q_X^* P_X^*$ such that $L \sigma = R$. Let $A$ be a finite automaton
accepting $L$, with vertices numbered $1, \dots, n$. Suppose without
loss of generality that the edges in $A$ are labelled by single
letters from $Q_X \cup P_X$. For each $i$ let $Q_i$ be the set of
all words in $Q_X^*$ which label paths from the initial vertex to
vertex $i$. Similarly, let $P_i$ be the set of all words in $P_X^*$
which label words from vertex $i$ to a terminal vertex.

Now if $w \in Q_i P_i$ then $w = uv$ where $u \in Q_X^*$ labels a
path from the initial vertex to vertex $i$, and $v \in P_X^*$ labels a
path from vertex $i$ to a terminal vertex. Hence $uv = w$ labels a path
from the initial vertex to a terminal vertex, and so $w \in L$.
Conversely, if $w \in L \subseteq Q_X^* P_X^*$ then $w$ admits a
factorisation $w = uv$ where $u \in Q_X^*$ and $v \in P_X^*$. Since
the edge labels in $A$ are single letters, an accepting path for $w$
must consist of a path from the initial vertex to some vertex $i$ labelled
$u$, followed by a path from $i$ to a terminal vertex labelled $v$.
It follows that $u \in Q_i$ and $v \in P_i$, so that $w \in Q_i P_i$.
Thus we have
$$L \ = \ \bigcup_{i=1}^n Q_i P_i$$
and so
$$R \ = \ L \sigma \ = \ \left( \bigcup_{i=1}^n Q_i P_i\right) \sigma \ = \ \bigcup_{i=1}^n \left( Q_i P_i\right) \sigma$$
as required.
\end{proof}

For the next proposition, we shall need some notation.
For a word $q = q_{x_1} q_{x_2} \dots q_{x_n} \in Q_X^*$, we let $q'
= p_{x_n} \dots p_{x_2} p_{x_1} \in P_X^*$. Similarly for a word $p
= p_{x_1} p_{x_2} \dots p_{x_n} \in Q_X^*$, we let $p' = q_{x_n}
\dots q_{x_2} q_{x_1} \in Q_X^*$. Note that $p'' = p$ and $q'' = q$.
Note also that $p' \sigma$ is the unique right inverse of $p \sigma$, and
$q' \sigma$ is the unique left inverse of $q \sigma$.

\begin{proposition}\label{factorise}
Let $u \in \Sigma_X^*$, and let $q \in Q_X^*$ and $p \in P_X^*$.
Then $u \sigma = (qp) \sigma$ if and only if there exists a
factorisation $u = u_1 u_2$ such that $(q'u_1) \sigma = 1 = (u_2 p')
\sigma$.
\end{proposition}
\begin{proof}
Suppose first that $u \sigma = (qp) \sigma$. Let $\Lambda$ be the
monadic rewriting system defined in the proof of Theorem~\ref{q*p*}.
Then $u$ is reduced by $\Lambda$ to $qp$. Notice that the only rules
in $\Lambda$ which can be applied to words not representing zero
remove factors representing the identity; it follows easily that $u$
admits a factorisation $u = u_1 u_2$ where $u_1 \sigma = q \sigma$
and $u_2 \sigma = p \sigma$. Now we have
$$(q' u_1) \sigma = (q' \sigma) (u_1 \sigma) = (q' \sigma) (q \sigma) = 1$$ and
symmetrically $(u_2 p') \sigma = 1$ as required.

Conversely, $q \sigma$ is the unique right inverse of $q' \sigma$,
so if $(q' u_1) \sigma =(q' \sigma) (u_1 \sigma) = 1$ then we must
have $u_1 \sigma = q \sigma$. Similarly, if $(u_2 p') \sigma = 1$ then $u_2
\sigma = p \sigma$, and so we deduce that $u \sigma = (u_1 u_2)
\sigma = (qp) \sigma$ as required.
\end{proof}

We are now ready to prove our main theorem about $M$-automata with
rational target sets where $M$ is a polycyclic monoid.

\begin{theorem}\label{polyzero}
Suppose $L \in F_{Rat}(P^0(X))$. Then $L$ is a finite union of
languages, each of which is the concatenation of one or two
languages in $F_1(P^0(X))$.
\end{theorem}

\begin{proof}
Let $M=P^0 (X)$ and let $A$ be an $M$-automaton with rational target
set $R$ accepting the language $L$. By Corollary~\ref{splitcor}
there exists an integer $n$ and regular languages $Q_1, \dots, Q_n
\subseteq Q_X^*$ and $P_1, \dots P_n \subseteq P_X^*$ such that
$$R = R_0 \cup \bigcup_{i=1}^n (Q_i P_i) \sigma.$$
where either $R_0 = \emptyset$ or $R_0 = \{ 0 \}$ depending on
whether $0 \in R$. For $1 \leq i \leq n$, we let $R_i = (Q_i P_i)
\sigma$. It follows easily that we can write
$$L = L_0 \cup L_1 \cup \dots \cup L_n$$
where each $L_i$ is accepted by a $M$-automaton with target set
$R_i$. Clearly it suffices to show that each $L_i$ is a finite union
of languages, each of which is the concatenation of at most two languages
in $F_1(M)$.

We begin with $L_0$. Let $Z = \lbrace u \in \Sigma_X^* \mid u \sigma
= 0 \rbrace$ and \mbox{$W = \lbrace w \in \Sigma_X^* \mid w \sigma =
1 \rbrace$}. It is easily seen (for example, by considering the
rewriting system $\Lambda$ from the proof of Theorem~\ref{q*p*})
that $u \in Z$ if and only if either $u$ contains the letter $z$, or
$u$ factorizes as $u_1 p_x u_2 q_y u_3$ where $x, y \in X$, $x \neq
y$ and $u_1,u_2,u_3 \in \Sigma_X^*$ are such that $u_2$ represents
the identity, that is, such that $u_2 \in W$. Thus,
$$Z \ = \ \Sigma_X^* \ \lbrace z \rbrace \ \Sigma_X^* \ \cup \ \bigcup_{x, y \in X, x \neq y} \Sigma_X^* \ \lbrace p_x  \rbrace \ W \ \lbrace q_y \rbrace \ \Sigma_X^*.$$
From this expression it is a routine matter to show that $Z$ is a
rational transduction of $W.$ By Proposition~\ref{transduction},
$L_0$ is a rational transduction of the language $Z$. Since the
class of rational transductions is closed under composition, it
follows that $L$ is a rational transduction of $W,$ and hence by
Proposition~\ref{transduction} that $L_0 \in F_1(M)$, as required.

We now turn our attention to the languages $L_i$ for $i \geq 1$.
Recall that $L_i$ is accepted by a $M$-automaton with target set
$R_i = (Q_i P_i) \sigma$. Let
$$P_i' \ = \ \lbrace (p', \epsilon) \mid p \in P_i \rbrace \ \subseteq \ Q_X^* \times \Sigma^*$$
and similarly
$$Q_i' \ = \ \lbrace (q', \epsilon) \mid q \in Q_i \rbrace \ \subseteq \ P_X^* \times \Sigma^*.$$
It is readily verified that $P_i'$ and $Q_i'$ are rational subsets
of $\Sigma_X^* \times \Sigma^*$; let $A_P$ and $A_Q$ be finite
automata accepting $P_i'$ and $Q_i'$ respectively, and assume
without loss of generality that the first component of every edge
label is either a single letter in $\Sigma_X$ or the empty word
$\epsilon$.

By Proposition~\ref{transduction} there is a rational transduction $\rho
\subseteq \Sigma_X^* \times \Sigma^*$ such that $w \in L_i$ if and
only if $(u, w) \in \rho$ for some $u \in \Sigma_X^*$ such that $u
\sigma \in R_i$. Let $A$ be an automaton recognizing $\rho$, again
with the property that the first component of every edge label is
either a single letter in $\Sigma_X$ or the empty word $\epsilon$.
We construct a new automaton $B$ with
\begin{itemize}
\item vertex set the disjoint union of the state sets of $A_Q$, $A$, and $A_P$;
\item all the edges of $A_Q$, $A$ and $A_P$;
\item initial vertex the initial vertex of $A_Q$;
\item terminal vertices the terminal vertices of $A_P$;
\item an extra edge, labelled $(\epsilon, \epsilon)$, from each terminal vertex of $A_Q$
to the initial vertex of $A$; and
\item an extra edge labelled $(\epsilon, \epsilon)$, from each terminal vertex of $A$ to
the initial vertex of $A_P$.
\end{itemize}
It is immediate that $B$ recognizes the relation
$$\tau \ = \ Q_i' \rho P_i' \ = \ \{ (q' x p', w) \mid q \in Q_i, p \in P_i, (x, w) \in \rho \} \ \subseteq \ \Sigma_X^* \times \Sigma^*$$
and again has the property that the first component of every edge
label is either a single letter or the empty word.

Let $Q$ be the vertex set of $A$, viewed as a subset of the vertex set
of $B$. For each vertex $y \in Q$, we let $K_y$ be the language of
all words $w$ such that $(u, w)$ labels a path in $B$ from the
initial vertex of $B$ to $y$ for some $u$ with $u \sigma = 1$. By
considering $B$ as an transducer but with terminal vertex $y$, we see
that $K_y$ is a rational transduction of the word problem of $P(X)$,
and hence by Proposition~\ref{transduction} lies in the class
$F_1(P(X))$.

Dually, we let $L_y$ be the language of all words $w$ such that $(u,
w)$ labels a path in $B$ from $y$ to a terminal vertex for some $u$
with $u \sigma = 1$. This time by considering $B$ as a transducer
but with initial vertex $y$, we see that $L_y$ is also a rational
transduction of the word problem of $P(X)$, and hence also lies in
$F_1(P(X))$.

We claim that
$$L_i = \bigcup_{y \in Q} K_y L_y,$$
which will clearly suffice to complete the proof.

Suppose first that $w \in L_i$. Then there exists a word $u \in
\Sigma_X^*$ such that $u \sigma \in R_i$ and that $(u, w) \in
\rho$. Since $R_i = (Q_i P_i) \sigma$ we have $u \sigma = (qp)
\sigma$ for some $q \in Q_i$ and $p \in P_i$.
 Note that $(q'up', w) \in \tau$ is accepted by $B$. By Proposition~\ref{factorise},
$u$ admits a factorization $u = u_1 u_2$ such that $(q'u_1) \sigma =
1$ and $(u_2 p') \sigma = 1$. Now in view of our assumption on the
edge labels of $B$, $w$ must admit a factorization $w = w_1 w_2$
such that $B$ has a path from the initial vertex to some vertex $y$
labelled $(q'u_1, w_1)$ and a path from $y$ to a terminal vertex
labelled $(u_2 p', w_2)$; moreover, the vertex $y$ can clearly be assumed to
lie in $Q$. Since $(q'u_1) \sigma = 1 = (u_2 p') \sigma$, it follows
that $w_1 \in K_y$ and $w_2 \in L_y$ so that $w = w_1 w_2 \in K_y
L_y$, as required.

Conversely, suppose $y \in Q$ and that $w = w_1 w_2$ where $w_1 \in
K_y$ and $w_2 \in L_y$. Then $B$ has a path from the initial vertex
to vertex $y$ labelled $(u_1, w_1)$ and a path from the vertex $y$ to
a terminal vertex labelled $(u_2, w_2)$ for some $u_1$ and $u_2$ with
$u_1 \sigma = u_2 \sigma = 1$. Since $y \in Q$, it follows from the
definition of $B$ that $u_1 = q' v_1$ and $u_2 = v_2 p'$ for some $q
\in Q_i$ and $p \in P_i$ and $v_1$ and $v_2$ such that $(v_1 v_2, w)
\in \rho$. But now $(q' v_1) \sigma = u_1 \sigma = 1$ and $(v_2
p') \sigma = u_2 \sigma = 1$, so we deduce by
Proposition~\ref{factorise} that $v_1 \sigma = q \sigma$ and $v_2
\sigma = p \sigma$. But then $(v_1 v_2) \sigma = (qp) \sigma \in
R_i \subseteq R$ and $(v_1 v_2, w) \in \rho$, from which it follows that $w
\in L_i$ as required.

Thus, we have written $L$ as a finite union of languages $L_i$ where
each $L_i$ either lies in $F_1(M)$ (in the case $i = 0$) or is a
finite union of concatenations of two languages in $F_1(M)$. This
completes the proof.
\end{proof}

In the case that $|X| \geq 2$, we have $P^0(X) = P(X)$ and $F_1(P(X))$ is the class of
context-free languages, which is closed under both finite union and
concatenation. Hence, we obtain the following easy consequence.
\begin{theorem}\label{contextfree}
If $|X| \geq 2$ then $F_{Rat}(P(X))$ is the class of context-free languages.
\end{theorem}

In the case $|X| = 1$, we have that $P^0(X)$ is isomorphic to the bicyclic
monoid $B = P(X)$ with a zero adjoined. Combining Theorem~\ref{polyzero}
with Proposition~\ref{M0M} and Theorem~\ref{G0G} we thus obtain.
\begin{corollary}
Every language in $F_{Rat}(B)$ is a finite union of languages, each
of which is the concatenation of one or two blind one-counter
languages.
\end{corollary}
Since the class $F_1(B)$ of partially blind one-counter languages is
not closed under concatenation, however, we cannot here conclude
that $F_{Rat}(B) = F_1(B)$. Indeed, the following result shows that
this is not the case.
\begin{theorem}
The language
$$\lbrace a^i b^i a^j b^j \mid i, j \geq 0 \rbrace$$
lies in $F_{Rat}(B)$ but not in $F_1(B)$.
\end{theorem}
\begin{proof}
Let $L = \{a^i b^i a^j b^j \mid i,j \geq 0 \}$. First, we claim that the
$B$-automaton with rational target set shown in Figure \ref{fig_L} accepts
the language $L$. Indeed, it is easily seen to accept exactly pairs of the
form
\[ (p^{i_{0}}q^{i_{1}}qpp^{i_{2}}q^{i_{3}},
a^{i_{0}}b^{i_{1}}a^{i_{2}}b^{i_{3}}) =
(p^{i_{0}}q^{i_{1}+1}p^{i_{2}+1}q^{i_{3}},a^{i_{0}}b^{i_{1}}a^{i_{2}}b^{i_{3}})\]
for $i_0, i_1, i_2, i_n \in \mathbb{N}$. A straightforward argument shows
that
$p^{i_{0}}q^{i_{1}+1}p^{i_{2}+1}q^{i_{3}} = qp$ if and only if
$i_0 = i_1$ and $i_2 = i_3$, which suffices to establish the claim and
proof that $L \in F_{Rat}(B)$.
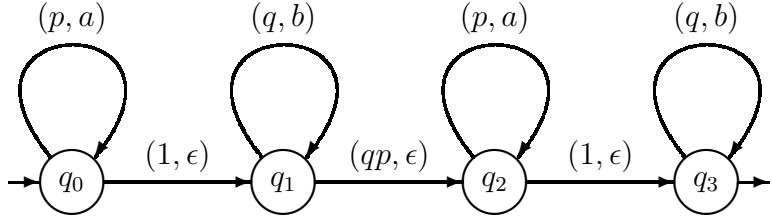
\begin{figure}
\begin{picture}(80,25)
\thicklines \setloopdiam{10} \Large
\letstate Q0=(10,5)    \drawinitialstate(Q0){$q_{0}$}
\letstate Q1=(30,5)       \drawstate(Q1){$q_{1}$}
\letstate Q2=(50,5)      \drawstate(Q2){$q_{2}$}
\letstate Q3=(70,5)        \drawfinalstate(Q3){$q_{3}$}
\drawloop[t](Q0){$(p, a)$} \drawloop[t](Q1){$(q, b)$}
\drawtrans(Q0,Q1){$(1, \epsilon)$} \drawloop[t](Q2){$(p,a)$}
\drawloop[t](Q3){$(q, b)$} \drawtrans(Q1,Q2){$(qp, \epsilon)$}
\drawtrans(Q2,Q3){$(1,\epsilon)$}
\end{picture}
\caption{A rational $B$-automaton with target set $\{qp\}$,
 accepting the language $\{a^i b^i a^j b^j \mid i,j \geq 0 \}$.}
\label{fig_L}
\end{figure}

Assume now for a contradiction that $L \in F_{1}(B)$. Then there exists
a $B$-automaton $A$ accepting $L$, with $N$ vertices say. For $i \geq 0$
let $\pi_{i}$ be an accepting path for $a^i b^i a^i b^i$. Suppose without
loss of generality that the right-hand sides of edge labels in $A$ are all
$a$, $b$ or $\epsilon$. Then we can write $\pi_{i} =
\alpha_{i}\beta_{i}\gamma_{i}\delta_{i}$ and where $\alpha_{i}$ has label
$( s_{i}, a^i)$,
$\beta_{i}$ has label $(t_{i},b^i)$, $\gamma_{i}$ has label $(u_{i}, a^i)$
and $\delta_{i}$ has label $(v_{i}, b^i)$ for some $s_i, t_i, u_i, v_i \in B$.

The proof will proceed by considering loops (that is, closed paths) in the
automaton $A$; we begin by introducing some terminology to describe particular
types of loops. A loop with label $(q^k p^j, x)$ is called
an \emph{increment loop} if $j>k$, a \emph{stable loop} if $k=j$ and a
\emph{decrement loop} if $k>j$. We call the loop an \emph{epsilon loop}
if $x=\epsilon$ and a \emph{non-epsilon} loop otherwise. A path which
does not traverse any loops is called a \emph{simple path}.

First, notice that since there are only finitely many simple paths, there
exists a constant $K$ such that every simple path in $A$ has label of the
form $(q^g p^h, x)$ with $g + h < K$.

Now let us consider paths of the form $\alpha_{i}$. We claim that for all
but at most $KN$ values of $i$, the path $\alpha_{i}$ contains a non-epsilon increment
loop. For all $i \geq N$, we can write $\alpha_i = \alpha_i^{(1)} \alpha_i^{(2)}$
where $\alpha_i^{(1)}$ has label $(s_i^{(1)}, a^{i-N})$ and
$\alpha_i^{(2)}$ has label $(s_i^{(2)}, a^N)$.

First note that the only elements of $B$ which generate a right
ideal [left ideal] including the identity element, are those of the
form $p^k$ [respectively $q^k$] for some $k \geq 0$. Thus, we must
have that both $s_i$ and $s_i^{(1)}$ are powers of $p$, and that
$v_i$ is a power of $q$. In particular, we can let $f_{i} \geq 0$ be
such that $s_i^{(1)} = p^{f_{i}}$.

First suppose $i$ is such that $\alpha_i^{(1)}$ does not traverse an
increment loop. Let $\alpha_{i}'$ be the path obtained from
$\alpha_i^{(1)}$ by removing all loops, and suppose $\alpha_{i}'$
has label $(q^g p^h, a^l)$. Since none of the loops removed were
increment loops, it follows easily that
$$f_i \ \leq \ h - g \ \leq \ h+g \ \leq \ K.$$
Suppose now for a contradiction than more than $KN$ values of $i
\geq N$ are such that $\alpha_i'$ contains no increment loop. Then
by the pigeonhole principle, there exist $i \neq j$ with $i \geq N$
and $j \geq N$ such that $f_i = f_j$ and the paths $\alpha_i^{(1)}$
and $\alpha_j^{(1)}$ end at the same state. But now the composition
$\alpha_i^{(1)} \alpha_j^{(2)} \beta_j \gamma_j \delta_j$ is an
accepting path with label
\begin{align*}
(s_i^{(1)} s_j^{(2)} t_j u_j v_j, a^{i-N} a^N b^j a^j b^j) &= (p^{f_i^{(1)}} s_j^{(2)} t_j u_j v_j, a^i b^j a^j b^j) \\
&= (s_j^{(1)} s_j^{(2)} t_j u_j v_j, a^i b^j a^j b^j) \\
&= (s_j t_j u_j v_j, a^i b^j a^j b^j) \\
&= (1, a^i b^j a^j b^j)
\end{align*}
so that $a^i b^j a^j b^j$ is accepted by $A$, giving a contradiction.
Thus, we have established that for all but $KN$ values of $i \geq N$, the path
$\alpha_{i}^{(1)}$ must traverse an increment loop. Hence, for all but
$KN + N = (K+1)N$ values of $i \geq 0$, the path $\alpha_{i}^{(1)}$ must
traverse an increment loop.

Now let $i$ be such that $\alpha_{i}^{(1)}$ traverses an increment
loop and suppose for a contradiction that $\alpha_i$ does not
traverse a non-epsilon increment loop. Consider the path
$\alpha_{i}^{(2)}$. Clearly, since this path has label with
right-hand-side $a^N$, and the right-hand-sides of edge labels in
the automaton are single letters or $\epsilon$, this path must
traverse a non-epsilon loop. Since $\alpha_{i}$ does not traverse a
non-epsilon increment loop,  $\alpha_{i}^{(2)}$ must traverse a
non-epsilon stable or decrement loop, say with  label $(q^g p^h,
a^{k})$ where $0 \leq h \leq g$ and $0 < k$. We also know that
$\alpha_{i}^{(1)}$ traverses an epsilon increment loop, say with
label $(q^x p^y, \epsilon)$ where $0 \leq x < y$. Clearly, by
traversing the latter loop an additional $(g-h)$ times and the
former loop an additional $(y-x)$ times, we obtaining an accepting
path for the word $a^{i + (y-x)k} b^i a^i b^i$, which gives the
required contradiction.

Thus, we have shown that for all but at most $(K+1)N$ values of $i$, the path
$\alpha_i$ traverses a non-epsilon increment loop.
A left-right symmetric argument can be used to establish firstly that
each $v_{i} = q^{g_{i}}$ for some $g_{i} \geq 0$, and then that for
$i$ sufficiently large, $\delta_{i}$ must traverse a non-epsilon decrement
loop. Thus, for all but at most $2(K+1)N$ values of $i$, the paths
$\alpha_i$ and $\delta_i$ traverse respectively a non-epsilon increment
loop and a non-epsilon decrement loop.

Now choose $i$ such that this holds, and let $(q^j p^k, a^m)$ label
a non-epsilon increment loop in $\alpha_{i}$ and let $(q^{j'} p^{k'}, b^{m'})$
label a non-epsilon decrement loop in $\delta_{i}$ where
$k > j$, $k' > j'$ and $m,m' >0$.
Let $\pi'_i$ be the path obtained from $\pi_i$ by traversing the
given increment loop an additional $j'-k'$ times, and the given decrement
loop an additional $k-j$ times. Then $\pi_i$ has label of the form
$$\left( \ t (q^j p^k)^{(j'-k')+1} u (q^{j'} p^{k'})^{(k-j)+1} v, \ a^{i+m(j'-k')} b^i a^i b^{i + m'(k-j)} \ \right)$$
where $t$, $u$ and $v$ are such that $\pi$ has label
$$\left( \ t q^j p^k u q^{j'} p^{k'} v, \ a^i b^i a^i b^i \ \right)$$
so that in particular $t q^j p^k u q^{j'} p^{k'} v = 1$.
Now by our argument above regarding right and left ideals, the element
$tq^j \in B$ must be a power of $p$, while $q^{j'}v \in B$ must be a power of $q$.
Noting that powers of $p$ commute with each other, and powers of $q$ commute with
each other, we get
\begin{align*}
t (q^j p^k)^{(j'-k')+1} u (q^{j'} p^{k'})^{(k-j)+1} v \
&= \ t q^j p^{(k-j)(j'-k')} p^k u q^{j'} q^{(k-j)(j'-k')}p^{k'}v \\
&= \ p^{(k-j)(j'-k')}t q^j p^k u q^{j'} p^{k'}v q^{(k-j)(j'-k')} \\
&= \ p^{(k-j)(j'-k')}1q^{(k-j)(j'-k')} \\
&= \ 1.
\end{align*}
Therefore $\pi_{i}'$ is an accepting path. Thus, the automaton
accepts the word \[a^{i+m(j'-k')} b^i a^i b^{i + m'(k-j)}\] which is
not in the language $L$, giving the required contradiction. This
completes the proof that $L \not\in F_{1}(B)$.
\end{proof}

It is possible, however, to describe concatenations of partially
blind one-counter languages using partially blind two-counter
automata. Indeed more generally we have the following proposition.
\begin{proposition}\label{product}
Let $M_1$ and $M_2$ be monoids and $L_1$ and $L_2$ languages over
the same alphabet. If $L_1 \in F_1(M_1)$ and $L_2 \in F_1(M_2)$ then
$L_1 L_2 \in F_1(M_1 \times M_2)$.
\end{proposition}
\begin{proof}
By Proposition~\ref{transduction} for $i = 1, 2$ there are alphabets
$\Omega_{i}$, morphisms $\omega_{i}:\Omega_i^* \to M_{i}$ and rational
transductions $\rho_i \subseteq \Omega_i^* \times \Sigma^*$ such
that $L_i = \{1\}\omega_{i}^{-1}\rho_{i}$. Assume without loss of
generality that $\Omega_1$ and $\Omega_2$ are disjoint, and let
$\Omega = \Omega_1 \cup \Omega_2$. Then there is a natural morphism
$\omega : \Omega^* \to M_1 \times M_2$ extending
$\omega_{1},\omega_{2}$. Now let $\rho$ be the product of $\rho_{1}$
and $\rho_{2}$:
$$\rho \ = \ \lbrace (u_1 u_2, w_1 w_2, \mid (u_1, w_1) \in \rho_1, (u_2, w_2 \in \rho_2) \rbrace \ \subseteq \ \Omega^* \times \Sigma^*.$$
Then $\rho$ is a rational transduction from $\Omega^*$ to $\Sigma^*$.
Clearly, if $u_1 \in \Omega_1^*$ and $u_2 \in \Omega_2^*$ then $u_1 u_2$
represents the identity element in $M_1 \times M_2$ if and only if
$u_1$ and $u_2$ represent the identity elements in $M_1$ and $M_2$ respectively.
It follows that $w$ is in the image under $\rho$ of the identity language of
$M_1 \times M_2$ if and only if $w = w_1 w_2$ where $w_1 \in L_1$ and
$w_2 \in L_2$, so that $w \in L_1 L_2$.
Thus, $L_1 L_2$ is a rational transduction of the identity language of
$M_1 \times M_2$, so applying Proposition~\ref{transduction} again
we see that $L_1 L_2 \in F_1(M_1 \times M_2)$ as required.
\end{proof}

Since classes of the form $F_1(M)$ are closed under union,
Proposition~\ref{bicyclic}, Theorem~\ref{polyzero}
and Proposition~\ref{product} combine to give the following inclusion.
\begin{corollary}
$$F_{Rat}(B) \subseteq F_1(B^2).$$
\end{corollary}

\section*{Acknowledgements}

The research of the second author was supported by an RCUK Academic
Fellowship.

\bibliographystyle{plain}

\end{document}